\documentstyle[12pt]{article}

\begin{document}

\newcommand{\fig}[4]{
        \begin{figure}[htbp]
        \setlength{\epsfysize}{#2}
        \centerline{\epsfbox{#4}}
        \caption{#3} \label{#1}
        \end{figure}
        }

\addtolength{\parskip}{1ex}

\def\a{\alpha}
\def\b{\beta}
\def\c{\chi}
\def\d{\delta}
\def\D{\Delta}
\def\e{\epsilon}
\def\f{\phi}
\def\F{\Phi}
\def\g{\gamma}
\def\G{\Gamma}
\def\k{\kappa}
\def\K{\Kappa}
\def\z{\zeta}
\def\th{\theta}
\def\Th{\Theta}
\def\l{\lambda}
\def\la{\lambda}
\def\m{\mu}
\def\n{\nu}
\def\p{\pi}
\def\P{\Pi}
\def\r{\rho}
\def\R{\Rho}
\def\s{\sigma}
\def\S{\Sigma}
\def\t{\tau}
\def\om{\omega}
\def\Om{\Omega}
\def\smallo{{\rm o}}
\def\bigo{{\rm O}}
\def\to{\rightarrow}
\def\E{{\bf Exp}}
\def\ex{{\bf Exp}}
\def\cd{{\cal D}}
\def\rme{{\rm e}}
\def\hf{{1\over2}}
\def\R{{\bf  R}}
\def\cala{{\cal A}}
\def\cale{{\cal E}}
\def\calz{{\cal Z}}
\def\cald{{\cal D}}
\def\calp{{\cal P}}
\def\Fscr{{\cal F}}
\def\cc{{\cal C}}
\def\calc{{\cal C}}
\def\calh{{\cal H}}
\def\calv{{\cal V}}
\def\bk{\backslash}

\def\out{{\rm Out}}
\def\temp{{\rm Temp}}
\def\overused{{\rm Overused}}
\def\big{{\rm Big}}
\def\moderate{{\rm Moderate}}
\def\swappable{{\rm Swappable}}
\def\candidate{{\rm Candidate}}
\def\bad{{\rm Bad}}
\def\crit{{\rm Crit}}
\def\col{{\rm Col}}
\def\dist{{\rm dist}}
\def\mod{{\rm mod }}
\def\supp{{\rm supp}}

\newcommand{\Exp}{\mbox{\bf Exp}}
\newcommand{\var}{\mbox{\bf Var}}
\newcommand{\pr}{\mbox{\bf Pr}}

\newtheorem{lemma}{Lemma}
\newtheorem{theorem}[lemma]{Theorem}
\newtheorem{corollary}[lemma]{Corollary}
\newtheorem{claim}[lemma]{Claim}
\newtheorem{fact}[lemma]{Fact}
\newtheorem{proposition}[lemma]{Proposition}
\newtheorem{prob}[lemma]{Problem}
\newtheorem{question}[lemma]{Question}
\newtheorem{conjecture}[lemma]{Conjecture}
\newtheorem{example}[lemma]{Example}
\newenvironment{proof}{{\bf Proof.}}{\hfill\rule{2mm}{2mm}}
\newenvironment{definition}{\noindent{\bf Definition.}}{}
\newtheorem{remarka}[lemma]{Remark}
\newenvironment{remark}{\begin{remarka}\rm}{\end{remarka}}

\newcommand{\limninf}{\lim_{n \rightarrow \infty}}
\newcommand{\proofstart}{{\bf Proof\hspace{2em}}}
\newcommand{\tset}{\mbox{$\cal T$}}
\newcommand{\proofend}{\hspace*{\fill}\mbox{$\Box$}}
\newcommand{\bfm}[1]{\mbox{\boldmath $#1$}}
\newcommand{\reals}{\mbox{\bfm{R}}}
\newcommand{\expect}{\mbox{\bf Exp}}
\newcommand{\he}{\hat{\e}}
\newcommand{\card}[1]{\mbox{$|#1|$}}
\newcommand{\rup}[1]{\mbox{$\lceil{ #1}\rceil$}}
\newcommand{\rdn}[1]{\mbox{$\lfloor{ #1}\rfloor$}}
\newcommand{\ov}[1]{\mbox{$\overline{ #1}$}}
\newcommand{\inv}[1]{\mbox{$1\over #1 $}}

\def\calp{{\cal P}}
\newcommand{\csp}{{CSP_{n,p}(\calp)}}
\newcommand{\cspc}{{CSP_{n,p=c/n^{k-1}}(\calp)}}
\newcommand{\cspp}{{CSP(\calp)}}
\newcommand{\csph}{{\widehat{\CSP}_{n,p}({\cal P})}}
\def\CSP {{CSP}}
\def\Dom {{\rm Dom}}
\def\supp {{\bf supp}}
\def\Prob {{\bf Pr}}
\def\Ex {{\bf Ex}}
\def\Var {{\bf Var}}
\def\thr {{\rm thr }}

\title{{\Large \bf Sharp thresholds for constraint satisfaction problems and homomorphisms}}

\author{
Hamed Hatami {\em and} Michael Molloy\\
Department of Computer Science\\
University of Toronto\\
Toronto, Canada\\
\{hamed,molloy\}@cs.toronto.edu
}
\maketitle

\begin{abstract}
We determine under which conditions certain natural models of random constraint satisfaction problems
have sharp thresholds of satisfiability.  These models  include graph and hypergraph
homomorphism, the $(d,k,t)$-model, and binary constraint satisfaction problems with domain size 3.
\end{abstract}
\thispagestyle{empty}
\newpage
\setcounter{page}{1}

\section{Introduction}
Random 3-SAT and its generalizations have been studied intensively for the past
decade or so (see eg.\ \cite{abm, ac6,cs,ef,jsv,nae,cr,dm2,fla,xl1,cm}).
One of the most interesting things about these models,
and arguably the main reason that most people study them, is that many of them
exhibit what is called a {\em sharp threshold of
satisfiability}\footnote{Defined formally below.}, a critical clause-density
at which the random problem suddenly moves from being almost surely\footnotemark[1] satisfiable
to almost surely unsatisfiable.    Most
of the work on these problems is, at least implicitly, an attempt
to determine the precise locations of their thresholds.  At this point, these locations
are known only for a handful of the problems, such as \cite{nae,cr,dm2,fla,xl1,cm}.
Just proving the
existence of a sharp threshold for random 3-SAT was considered a major breakthrough
by Friedgut\cite{ef}.  The vast majority of these generalizations appear to
have sharp thresholds, but there are exceptions which are said to have coarse
thresholds\footnotemark[1].

The ultimate goal of the present line of enquiry is to determine
precisely which of these models have sharp thresholds, but this appears to be quite difficult;
in Section \ref{s3} we show that it is at least as difficult as determining
the location of the threshold for 3-colourability, something that has
been sought after for more than 50 years (see eg.\ \cite{er,mmsurv}).
A more fundamental goal is to obtain a better understanding of what
can cause some
problems to have coarse thresholds rather than sharp ones.

Molloy\cite{mm1} and independently Creignou and Daude\cite{cd}
introduced a wide family of models for random constraint satisfaction
problems which includes 3-SAT and many of its generalizations.
This permits us to study them under a common
umbrella, rather than one-at-a-time. Molloy determined precisely which models from
this family have any threshold at all (\cite{cd} provides the
same result for those models with domain\footnotemark[1] size 2).
But he left open
the much more important question of which models have sharp thresholds.
In this paper, we begin to address this question.  We answer it for two of the most
natural subfamilies - the so-called $(d,k,t)$-family\footnotemark[1]
(Theorem \ref{tdkt}),
and the family of graph and hypergraph homomorphism problems (Theorem \ref{thomo}).
We also shed light
on the more fundamental problem by determining the only properties that can
cause a coarse threshold in binary constraint satisfaction problems
with domain size 3.

The standard example of a problem with a coarse threshold is
2-colourability. Here, there is a coarse threshold precisely
because unsatisfiability (i.e. non-2-colourability) can be caused
only by the presence of odd cycles.  Roughly speaking, Friedgut's
theorem\cite{ef} implies that a problem exhibits a coarse
threshold iff unsatisfiability is {\em approximately} equivalent
to having one of a set of unicyclic\footnotemark[1] subproblems.
It is not hard to see that if there are unsatisfiable unicyclic
instances of a problem then that problem exhibits a coarse
threshold (or exhibits no threshold at all). This makes it quite
natural to pose the following rule-of-thumb:

{\bf Hypothesis A:} {\em If a random model from the family in \cite{mm1} is such that:
(a) it exhibits a threshold, and (b) every unicyclic instance is satisfiable,
then that threshold is sharp.}

However, reality is not that simple. \cite{mm1} presents a counterexample
to Hypothesis A; others are presented in this paper.  Nevertheless,
the hypothesis holds for certain subfamilies of models. Creignou and Daude\cite{cd}
conjectured that Hypothesis A holds for problems with domain-size
two;  this was proven by Istrate\cite{gi} and independently
Creignou and Daude\cite{cd2} proved it for the case where the model is
{\em symmetric}.  Theorems \ref{tdkt} and \ref{thomo} in this paper
show that Hypothesis A holds for the
$(d,k,t)$-models and for homomorphism problems.

In general, coarse thresholds can be caused by
much more subtle and insidious reasons than unsatisfiable unicyclic instances.
In this paper we begin to understand some of these reasons by focusing on the
case where the constraint size is two and the domain size is three (a natural
next step after the well-understood domain-size-two case).  In this paper, we identify
a particular subtle cause, and show that this and unsatisfiable unicyclic instances
are the only things that can cause a coarse threshold (Theorem \ref{t23}).  If we
permit either greater domain sizes or greater constraint sizes then this
is no longer true - there are other possible causes.

\subsection{The random models}
In our setting, the variables of a constraint satisfaction problem (CSP)
all have the same domain of permissable
values, $\{1,...,d\}$, and all constraints will have size $k$,
for some fixed integers $d,k$.
Given a $k$-tuple of variables, $(x_1,...,x_k)$,
a {\em restriction} on $(x_1,...,x_k)$ is a  $k$-tuple of values
$R=(\d_1,...\d_k)$ where each $1\leq\d_i\leq d$.
For each $k$-tuple $(x_1,...,x_k)$, the set of restrictions
on that $k$-tuple is called a {\em constraint}.
The {\em empty constraint} is the constraint which contains
no restrictions. We say that an
assignment of values to the variables of a constraint $C$
{\em satisfies} $C$ if that assignment is not one of the restrictions
in $C$.   An assignment of values to all variables in a CSP {\em satisfies}
that CSP if every constraint is simultaneously satisfied.
A CSP is {\em satisfiable} if it has such a satisfying assignment.

It will be convenient to consider a set of canonical variables
$X_1,...,X_k$ which are used only to describe the ``pattern" of a constraint.
These canonical variables are not variables of the actual CSP.
For any $d,k$ there are $d^k$ possible restrictions and $2^{d^k}$
possible constraints over the $k$ canonical variables. We denote this
set of constraints as $\calc^{d,k}$.
For our random model, one begins by specifying a particular probability
distribution, $\calp$ over $\calc^{d,k}$.  We use $\supp(\calp)$ to
denote the support of $\calp$; i.e. the set of constraints $C$ with
$\calp(C)>0$.
Different choices of $\calp$ give rise to different instances of the
model.

We now define our random models.  The ``$G_{n,M}$'' model,
where the number of constraints is fixed to be $M$, is the most common.
But in this paper,
it will be much more convenient to focus on the ``$G_{n,p}$''
model where each $k$-tuple of variables is chosen independently with probability
$p=c/n^{k-1}$ to receive a constraint.  The two models are, in most respects,
equivalent when $M=(c/k!)n$.  In particular, it is straightforward to show that one
exhibits a sharp threshold iff the other does.

{\bf The $\csp$ Model:} Specify $n,p$ and $\calp$ (typically $p=c/n^{k-1}$
for some constant $c$;
note that $\calp$ implicitly specifies $d,k$).
First choose a random $k$-uniform
hypergraph on $n$ variables where each of the ${n\choose k}$ potential
hyperedges is selected with probability $p$.
Next, for each hyperedge $e$, we choose a constraint on the
$k$ variables of $e$ as follows: we take a random permutation from
the $k$ variables onto $\{X_1,...,X_k\}$ and then we select a random
constraint  according to $\calp$ and map it onto
the $k$ variables.

A property holds {\em almost surely} (a.s) if the limit as
$n\rightarrow\infty$ of it holding is $1$.  We say that $\csp$ has
a {\em sharp threshold of satisfiability} if there is some
$c=c(n)>0$ such that for every $\e>0$, if $p=(1-\e)c/n^{k-1}$ then
$\csp$ is a.s. satisfiable and if $p=(1+\e)c/n^{k-1}$ then $\csp$
is a.s. unsatisfiable.  This is often abbreviated to just {\em
sharp threshold}.  We say that $\csp$ has a {\em coarse threshold}
if for all $c$ in some interval $c_1(n)<c<c_2(n)$,  $\csp$ is
neither a.s. satisfiable nor a.s. unsatisfiable. If $\csp$ has
neither a sharp nor a coarse threshold, then it is easy to see
that it must either be a.s. satisfiable for all $c>0$ or a.s.
unsatisfiable for all $c>0$.

Each $k$-tuple of vertices can have at most one constraint in $\csp$.  When applying
Friedgut's theorem, it will be convenient to relax this condition, and allow
$k$-tuples to possibly receive multiple constraints.  Thus up to
$k!\times|\supp(\calp)|$ constraints can appear
on a $k$-tuple of variables.

{\bf The $\csph$ Model:} Specify $n,p$ and $\calp$. For each of the $n(n-1)...(n-k+1)$
{\em ordered} $k$-tuples of variables and each constraint $C\in\supp(\calp)$,
we assign $C$ to the ordered $k$-tuple with probability $\calp(C)\times p/k!$.

Note that the expected total number of constraints is the same under each model.
Furthermore, it is easy to calculate that the probability of at least one
$k$-tuple receiving more than one constraint in $\csph$ is for $k\geq 3$, $o(1)$
and for $k=2$, an absolute constant $0<\a< 1$.  It follows that if a property
holds a.s. in $\csph$ then it holds a.s. in $\csp$.  As a corollary, we have:

\begin{lemma}\label{l1}
If $\csph$ has a sharp threshold then so does $\csp$.
The reverse is true for $k\geq3$.
\end{lemma}

So for the remainder of the paper, whenever we wish to prove that $\csp$ has a sharp
threshold, we will work in the $\csph$ model.

We often focus on the {\em constraint hypergraph} of a CSP; i.e. the hypergraph
whose vertices are the variables and whose edges are the tuples of variables
that have constraints.  A {\em tree-CSP} is a CSP whose constraint
hypergraph is a hypertree. A CSP is {\em unicyclic} if its constraint
hypergraph is unicylic; i.e. has exactly one cycle. (Hypertree and cycle are defined below).

We close this subsection with some hypergraph definitions.
A hypergraph consists of a set of vertices and a set of {\em hyperedges},
where each hyperedge is a collection of vertices.  If every hyperedge
has size exactly $k$ then the hypergraph is {\em $k$-uniform}.  In a {\em simple}
hypergraph, no vertex appears twice in any one hyperedge, and no two edges
are identical.  So, for example, the constraint hypergraph of $\csp$ is
simple, but the constraint hypergraph of $\csph$ may have multiple edges.
Neither model permits multiple copies of a vertex in a single edge,
but such edges are possible when we discuss hypergraph homomorphism problems.
The edge $(v,v,...,v)$ is called a {\em loop}.

A walk $P$ of size $r$ is a sequence of $r$ hyperedges and $r+1$
vertices $(v_0, e_1,v_1, e_2, v_2 \ldots, e_r, v_r)$ such that
$e_i$ contains both $v_{i-1}$ and $v_i$.  A walk  is a {\em path}
if the $v_i$ are distinct.  A walk is a {\em cycle} of size $r$ if
for $i=1,\ldots,r$ the $v_i$ and $e_i$ are distinct, and
$v_0=v_r$. The {\em distance} from a vertex $u$ to a vertex $v$ is
the minimum $r$ such that there exists a walk of length $r$,
$(v,e_1,v_1,\ldots,e_r,u)$; the distance of a vertex from itself
is defined to be $0$. The distance from a vertex $v$ to a set of
vertices is the minimum distance from $v$ to any vertex in the
set. A hypergraph is a {\em hypertree} if it has no cycles and it
is connected.

By {\em contracting}
two vertices $u$ and $v$ into a new vertex $w$, we mean (i)
adding a new vertex $w$ to the set of the vertices, (ii) replacing
$u$ and $v$ in every hyperedge by $w$, and (iii) removing $u$ and $v$.

\subsection{Two special families}\label{section:hom}
Perhaps the most natural choice for $\calp$ is the distribution
obtained by selecting each of the $d^k$ possible restrictions
independently with probability $1/d^k$.
 However, as noted in \cite{akk}, every such
choice of $\calp$ yields a model that is a.s. unsatisfiable for
any non-trivial choice of $p$.  So this is a rather uninteresting
family of models, particularly as far as the study of thresholds
goes.

The next most natural choice for $\calp$ is to fix $t$, the number of restrictions
per clause, and to make every constraint with exactly $t$ restrictions equally likely.
(Note that for $d=2,t=1$ this yields random $k$-SAT.)
This is often called the $(d,k,t)$-model and has received a great deal of study,
both from a theoretical perspective \cite{mit,ms} and from experimentalists
(see \cite{gmp} for a survey of many such studies).  In \cite{akk} it is shown
that when $t\geq d^{k-1}$,  this model is problematic in the same way as
the previously mentioned one, as it is a.s. unsatisfiable even for values
of $p=o(1/n^{k-1})$ (i.e. when the number of constraints is $o(n)$).
However, it was proven in \cite{gmp} that for every $1\leq t<d^{k-1}$, the
$(d,k,t)$-model does not have that problem.  One of the main contributions of
this paper is to show that in this case the model  exhibits a sharp threshold:

\begin{theorem} \label{tdkt} For every $d,k\geq2$ and every $1\leq t< d^{k-1}$,
the $(d,k,t)$-model has a sharp threshold.
\end{theorem}

From a different perspective, it is quite natural to consider the case where
every constraint is identical, i.e. $|\supp(\calp)|=1$.  It is not hard
to see that every such problem
is equivalent to a hypergraph homomorphism problem, as defined below:

For two $k$-uniform hypergraphs, $G,H$, a {\em homomorphism} from
$G$ to $H$ is a mapping $h$ from $V(G)$ to $V(H)$  such that for
each edge $(v_1,v_2,\ldots,v_k)$ of $G$,
$(h(v_1),h(v_2),\ldots,h(v_k))$ is an edge of $H$. We say that $G$
is {\em homomorphic} to $H$, if there exists such a homomorphism.
When $k=2$ and $H$ is the complete graph with no loops, we are
simply asking whether $G$ has a $d$-colouring. Homomorphisms are
an important generalization of graph colouring (see, eg.
\cite{hn}).  They are often also referred to as $H$-colourings
(eg. \cite{hn2,gkp}).

Suppose that $H$ is a fixed hypergraph, and $G$ is a random
hypergraph on $n$ vertices where each of the ${n\choose k}$
potential hyperedges is selected with probability $p$.  Set $d$ to
be equal to the number of vertices in $H$ and define a constraint
$C$ with domain size $d$ and constraint size $k$ by saying that
$C$ permits $X_1=\d_1,...,X_k=\d_k$ iff $(\d_1,...,\d_k)$ is a
hyperedge of $H$.  Treat each vertex of $G$ as a variable with
domain $\{1,..,d\}$ and assign $C$ to each hyperedge of $G$. We
call this the {\em $H$-homomorphism problem}.

Thus we have an instance of $\csp$ where $C$ is the only
constraint in $\supp(\calp)$ and furthermore $C$ is symmetric
under permutations of the canonical variables; in other words, all
constraints are identical even under permutations of variables. It
is easy to see that every such $\calp$ corresponds to a
homomorphism problem; just take $H$ to be the hypergraph where
$(\d_1,...,\d_k)$ is a hyperedge iff $C$ permits
$X_1=\d_1,...,X_k=\d_k$.  Note that here a hyperedge in $H$ may
contain multiple copies of a vertex.

Thus, these $H$-homomorphism problems are not only important
as a fundamental graph problem, but also because they form a very natural subclass of our
family of random CSP models.
In this paper, we prove that Hypothesis A holds for every connected $H$.

It is easy to see that if $H$ has a loop $(\d,\d,...,\d)$
then every hypergraph  is trivially homomorphic to $H$
(just map every vertex to $\d$); so the $H$-homomorphism problem has no
threshold at all.  The other trivial case is where $H$ has
no hyperedges at all and so no non-trivial hypergraph
has an $H$-homomorphism.

\begin{lemma} \label{cycle-hom}Suppose that $H$ is a nontrivial hypergraph with
no loops. We have the following:

\begin{enumerate}
\item For $k\geq 3$, every unicylic hypergraph is homomorphic to a
single hyperedge, and hence to $H$.
\item For $k=2$: if the triangle is homomorphic to $H$, then so is
every unicyclic graph; and the triangle is homomorphic to $H$ iff
$H$ contains a triangle.
\end{enumerate}
\end{lemma}
\begin{proof}
To prove part (1), let $(v_0, e_1,v_1, e_2, v_2 \ldots, e_r, v_0)$
be the unique cycle of the hypergraph, and let
$(w_0,\ldots,w_{k-1})$ be a single hyperedge. Define
$h(v_i)=w_{(i\ \mod \ 2)}$, for every $0 \le i \le r-2$ and
$h(v_{r-1})=w_{2}$. It is easy to see that one can extend $h$ to a
homomorphism from the unicyclic hypergraph to $H$.

Part (2) easily follows from the easy and well-known fact that
every cycle is homomorphic to the triangle, and the triangle is
not homomorphic to any cycle of size greater that $3$.
\end{proof}

From Lemma~\ref{cycle-hom} we conclude that proving that
Hypothesis A holds whenever $H$ is connected and undirected is
equivalent to proving:

\begin{theorem}\label{thomo} If $H$ is a connected undirected loopless hypergraph
with at least one edge,
then the $H$-homomorphism problem has a sharp threshold iff (a) $k\geq 3$
or (b) $k=2$ and $H$ contains a triangle.
\end{theorem}

We do not have a strong feeling
as to whether the ``connected'' condition is necessary here; we discuss
the possibility of extending Theorem \ref{thomo} to disconnected graphs in
Section \ref{sh}.

\subsection{Tools}
Our main tool is distilled from Friedgut's main theorem in
\cite{ef}.  Friedgut reported to us\cite{ef2} that his proof can
be adapted to the setting of this paper. To provide Friedgut's
theorem for CSP's in its full power instead of being restricted to
the unsatisfiability property, we consider, as Friedgut did, every
monotone property where a property $A$ is called \emph{monotone}
if it is preserved under constraint addition. A property $A$ on
CSP's is called \emph{monotone symmetric} if it is monotone and
invariant under CSP automorphisms. For a property $A$, $A_n$
denotes the restriction of $A$ on CSP's with exactly $n$
variables.  Roughly speaking, Friedgut's theorem says that for a
value of $p$ that is ``within'' the coarse threshold, there is a
constant sized instance $M$ such that $\tau<\Pr[M \subseteq
\widehat{\CSP}_{n,p}({\cal P})]<1-\tau$ for some constant $\tau$
which does not depend on $n$, and adding $M$ to our random CSP
boosts the probability of being in $A$ by at least $2\a>0$,
whereas adding a linear number of new random constraints only
boosts it by at most $\a$. First, we must formalize what we mean
by ``adding $M$''. Given two CSP's $M,F$ where $M$ has $r$
variables, and $F$ has at least $r$ variables, we define $F\oplus
M$ to be the CSP obtained by choosing a random $r$-tuple of
variables in $F$ and then adding $M$ on those $r$ variables. Now
we can state Friedgut's theorem formally:

\begin{theorem}
\label{appendix:toolcsp} Let $A=\{A_i\}$ be a series of monotone
symmetric properties in $\widehat{\CSP}_{n,p}({\cal P})$ with a
coarse threshold. There exist, $p=p(n)$, $\tau,\alpha,\epsilon>0$,
a CSP $M$ whose constraints are chosen from $\supp({\cal P})$ such
that for an infinite number of $n$:

\begin{itemize}
\item[(a)] $\alpha<\Pr[\widehat{\CSP}_{n,p}({\cal P}) \in
A]<1-3\alpha$.

\item[(b)] $\tau<\Pr[M \subseteq \widehat{\CSP}_{n,p}({\cal
P})]<1-\tau$.

\item[(c)] $\Pr[\widehat{\CSP}_{n,p}({\cal P}) \oplus M \in A]
> 1- \alpha$.

\item[(d)] $\Pr[\widehat{\CSP}_{n,p(1+\epsilon)}({\cal P}) \in A]
< 1- 2\alpha$.
\end{itemize}
\end{theorem}

When as in our setting $p(n)=c(n)/n^{k-1}$,
Theorem~\ref{appendix:toolcsp}(b) implies that $M$ is a unicycle
CSP. So we obtain the following corollary which is our main tool
in this paper.

\begin{corollary}
\label{toolcsp} For any $\calp$, if $\csph$ has a coarse
threshold of satisfiablity then there exist $p=p(n)$,
$\alpha,\epsilon>0$, and a unicyclic CSP $M$ on a constant number of variables whose
constraints are chosen from $\supp(\calp)$
such that:

\begin{itemize}
\item[(a)] $\alpha<\Pr[\widehat{\CSP}_{n,p}({\cal P}) \mbox{ is unsatisfiable}]<1-3\alpha$.
\item[(b)] $\Pr[\widehat{\CSP}_{n,p(1+\epsilon)}({\cal P}) \mbox{ is unsatisfiable}] < 1- 2\alpha$.
\item[(c)] $\Pr[\widehat{\CSP}_{n,p}({\cal P}) \oplus M \mbox{ is unsatisfiable}]> 1- \alpha$.
\end{itemize}
\end{corollary}


Our next tool proves some properties for local parts of a random CSP.

\begin{lemma}
\label{local} Suppose that $p<cn^{1-k}$ for some positive constant
$c$, and let $G$ be an instance of $\widehat{\CSP}_{n,p}({\cal
P})$. Choose a set $T$ of $t$ random variables. Then for every
$\epsilon>0$, and integer $r>0$ there exists an integer
$L(c,t,r,\epsilon)$ such that with the probability of at least
$1-\epsilon$:
\begin{itemize}
\item[(i):] No constraint of $G$ contains more than one variable of $T$.
\item[(ii):] $G$
induces a forest on the set of the variables that are of distance
at most $r$ from $T$. \item[(iii):] There are at most $L$
variables that are of distance at most $r$ from $T$.
\end{itemize}
\end{lemma}
\begin{proof}
Let $E_1$, $E_2$, and $E_3$ denote the events $(i)$, $(ii)$, and
$(iii)$ respectively. Trivially
\begin{equation}
\label{eq1} \Pr[E_1] \ge 1 - \sum_{i=2}^k n^{k-i}{t \choose i}k! p=1-o(1).
\end{equation}
The expected number of the cycles of size at most $2r$ which
contain at least one variable in $T$ is at most $t\sum_{i=2}^{2r}
n^{ik-i-1}p^{i}$. Thus
\begin{equation}
\label{eq2} \Pr[E_2] \ge 1-t\sum_{i=2}^{2r} n^{ik-i-1}p^{i} \ge
1-\frac{2tr(1+c)^{2r}}{n}=1-o(1).
\end{equation}

The expected number of the variables in a distance of at most $r$
from $T$ is at most $t\sum_{i=1}^{r}n^{ik-i}p^i$. So by
Chebychev's inequality, for sufficiently large $L$:
\begin{equation}
\label{eq3} \Pr[E_3] \ge 1-\frac{t\sum_{i=1}^{r}n^{ik-i}p^i}{L}
\ge 1-\frac{\epsilon}{2}.
\end{equation}
The lemma follows from (\ref{eq1}), (\ref{eq2}), and (\ref{eq3}).
\end{proof}

Our third tool is easily proven with a straightforward first
moment calculation and concentration argument (via eg. the second
moment method or Talagrand's inequality).

\begin{lemma}
\label{tree} Let $T$ be a tree-CSP whose constraints are in
$\supp({\cal P})$. There exists $z=o(n^{1-k})$ such that a.s.
$\widehat{\CSP}_{n,z}({\cal P})$ contains $T$ as a sub-CSP.
\end{lemma}

\section{Difficulty}\label{s3}

The ultimate goal of this research is to characterize all
distributions ${\cal P}$ for which $\csp$
exhibits a sharp threshold. However, the following example indicates
that this is very difficult, even for binary CSP's (the case where $k=2$).
In particular, it is at least as difficult as determining the location
of the 3-colourability threshold, a heavily pursued open problem.
(The existence of that threshold was proven in \cite{af};
see \cite{mmsurv} for a recent survey and see \cite{ac7,mcd} for the
best current bounds on its location.)

We set $d=5$ and $k=2$, and define  two constraints by listing
their  pairs of forbidden values:
\begin{eqnarray*}
C_1&=&\{(1,1),(2,2)\} \cup (\{1,2\}\times\{3,4,5\})\cup (\{3,4,5\}\times\{1,2\}),\\
C_2&=&\{(3,3),(4,4),(5,5)\} \cup (\{1,2\}\times\{3,4,5\})\cup (\{3,4,5\}\times\{1,2\}).
\end{eqnarray*}
Note that each constraint forces the endpoints of every edge to take values that are either
both in $\{1,2\}$ or both in $\{3,4,5\}$.  A $C_1$ constraint says that they
have to be different values if they are both in $\{1,2\}$.  A $C_2$ constraint says
that they have to be different values if they are both in $\{3,4,5\}$.

We let $C_1$ occur with probability $q$ and $C_2$ occur with
probability $1-q$ in ${\cal P}$.
Set $c(q)=(1-q)/q$.

\begin{fact}\label{f3c}
\begin{enumerate}
\item[(a)] If $G_{n,p=c(q)/n}$ is a.s. 3-colourable, then $\csp$ has a sharp threshold.
\item[(b)] If there is some $\e>0$ such that $G_{n,p=(c(q)-\e)/n}$
is a.s. not 3-colourable, then $\csp$ has a coarse threshold.
\end{enumerate}
\end{fact}

Thus, determining the type of the
threshold for all such models
$\CSP_{n,p}({\cal P})$ requires
the knowledge of for which values of $c$, $G(n,\frac{c}{n})$ is a.s.
$3$-colourable, and for which values it is a.s. not $3$-colourable.

\proofstart
Choose our $\csp$ by first taking $G_{n,p=c/n}$ and
then setting each edge to be $C_1$ with probability $q$ and $C_2$ otherwise.
Let $G_1,G_2$ be the subgraphs formed by the edges chosen to be $C_1,C_2$ respectively.
If $c<1$ then all components of $G_{n,p=c/n}$ are trees or unicycles
and the CSP is trivially satisfiable.  So we
can focus on the range $c>1$ and we let $T$ denote
the giant component of $G_{n,p=c/n}$.  Note that
the variables of $T$ must either all take values from $\{1,2\}$ or
all take values from $\{3,4,5\}$.

{\em Case 1: $c>\inv{q}$.} Then $G_1$ is equivalent to $G_{n,p=c_1/n}$
for some $c_1>1$ and it follows easily that a.s. $G_1$ contains
a giant component which is not 2-colourable. This giant component is a subgraph
of $T$ and so the variables of $T$ must all take values from $\{3,4,5\}$.
If follows that the CSP is satisfiable iff $G_2$ is 3-colourable. Note
that  $G_2$ is equivalent to $G_{n,p=c_2/n}$
for some $c_2>c(q)$.

{\em Case 2: $c<\inv{q}$.}
Then $G_1$ is equivalent to $G_{n,p=c_1/n}$
for some $c_1<1$ and  $G_2$ is equivalent to $G_{n,p=c_2/n}$
for some $c_2<c(q)$.  If $G_2$ is a.s. 3-colourable then the CSP
is a.s. satisfiable. If $G_2$ is a.s. not 3-colourable then
the CSP is satisfiable iff $T$ is 2-colourable; i.e., if
$G_2$ does not have an odd cycle lying within $T$.  It is
easy to see that this occurs with probability between $\z$ and $1-\z$
for some $\z>0$; i.e. that the CSP is neither a.s. satisfiable nor
a.s. unsatisfiable.

Fact \ref{f3c} now follows.  If $G_{n,p=c(q)/n}$ is a.s. 3-colourable, then $\csp$ has a sharp threshold
which lies somewhere above $\inv{q}$. If there is some $\e>0$ such that $G_{n,p=(c(q)-\e)/n}$
is a.s. not 3-colourable, then $\csp$ has a coarse threshold running from $\inv{q}-\d$ to $\inv{q}$
for some $\d>0$.
\proofend

\section{Homomorphisms}\label{sh}

In this section, we prove our theorem concerning $H$-homomorphisms.
Let $G^k_{n,p}$  denote the random $k$-uniform hypergraph on $n$ vertices
where each $k$-tuple is present as a hyperedge with probability $p$.

{\bf Proof of Theorem \ref{thomo}} We begin with the case
$k\geq3$. Let $H$ be some $k$-uniform hypergraph, and assume that
$H$ has a coarse threshold.  Let $M,p,\a,\e$ be as guaranteed by
Corollary~\ref{toolcsp}.  In this setting, $M$ is a unicyclic
$k$-uniform hypergraph, such that adding $M$ to $G^k_{n,p}$ boosts
the probability of not having a homomorphism to $H$ by at least
$2\a$.

Consider $G=G^k_{n,p} \oplus M$. Let $M^+$ be the subgraph of $G$
consisting of all hyperedges that contain at least one vertex of
$M$ (and, of course all vertices in those hyperedges); in other
words, $M^+$ is the subhypergraph induced by the vertices of $M$
and all their neighbours. Lemma~\ref{local} implies that there is
some constant $L$ such that with probability at least
$1-\frac{\a}{2}$: $M^+$ is unicyclic and has at most $L$ vertices,
and no hyperedge of $G^k_{n,p}$ contains more than one vertex of $M$.

Since $M$ is unicylic, by Lemma~\ref{cycle-hom} there exists a
homomorphism $h$ from $M$ to a single edge, say
$(v_1,\ldots,v_k)$. Let $h_i$ be the set of the vertices in $M$
that are mapped by $h$ to $v_i$. Obtain the hypergraph $G'$ from
$G$ by (i) removing all edges in $M$; (ii) contracting all of the
vertices in $h_i$ into one single new vertex $u_i$, for each $1
\le i \leq k$; (iii) adding the single hyperedge
$(u_1,\ldots,u_k)$.

Suppose that $h'$ is a homomorphism from $G'$ to $H$. Then a
mapping from the vertices of $G$ to the vertices of $H$ which maps
every vertex $v$ in $G-M$ to $h'(v)$, and every vertex in $h_i$
to $h'(u_i)$ is a homomorphism from $G$ to $H$.  Thus, if $G'$ is
homomorphic to $H$ then so is $G$.

Let $T$ be the hypertree defined as follows: $T$ has a hyperedge
$(t_1,...,t_k)$, and each $t_i$ lies in $L$ other hyperedges.
Only $t_1,...,t_k$ lie in more than one edge of $T$.  Thus, $T$
has $k+k(k-1)L$ vertices and $kL+1$ hyperedges. Note that in $G'$
the subgraph induced by all edges containing $\{u_1,\ldots,u_k\}$
form a subtree of $T$. It follows that $G^k_{n,p}\oplus T$ is at
least as likely to be non-homomorphic to $H$ as $G'$ is, so:
$$\Pr[G^k_{n,p} \oplus M \mbox{ is not homomorphic to $H$}]
\le \Pr[G^k_{n,p} \oplus T \mbox{ is not homomorphic to $H$} ] + \frac{\alpha}{2}.$$
By Lemma~\ref{tree}, increasing $p$ by an additional $\e p$ a.s. results
in the addition of a copy of $T$. Thus:
\[\Pr[G^k_{n,p} \oplus T \mbox{ is not homomorphic to $H$}] \le
\Pr[G^k_{n,p(1+\epsilon)}\mbox{ is not homomorphic to $H$}]\]
which yields a contradiction to Corollary~\ref{toolcsp}(b).

This proves the case where $k\geq 3$, so we now turn to the case
$k=2$.  If $H$ contains no triangle, then $K_3$ is not homomorphic
to $H$.  Thus, $K_3$ forms a unicyclic unsatisfiable CSP using the
$H$-colouring constraints and so we do not have a sharp threshold.
So we will focus on graphs $H$ that contain a triangle.  Our proof
follows along the same lines as the case $k\geq 3$, but is
complicated a bit since we can no longer assume that $M$ is
homomorphic to a single edge. We only highlight the differences.

Define $M^+$ to be the subgraph of $G=G^k_{n,p}\oplus M$ induced
by all vertices within distance $r=|V(H)|+|V(M)|+3$ of the unique
cycle of $M$.  By Lemma~\ref{local} there is some constant $L$
such that with probability at least $1-\frac{\a}{2}$: $M^+$ is
unicyclic and has at most $L$ vertices.

Define $U$ to be the set of vertices of $G$ that are of distance
exactly $r=|V(H)|+|V(M)|+3$ from the unique cycle of $M$. Consider
any vertex $u\in H$. By Lemma~\ref{distance} below, if $M^+$ is unicyclic
then there is a homomorphism from $M^+$ to $H$ such that all
vertices in $U$ are mapped to $u$.

Obtain the graph $G'$ from $G$ by (i) removing all
of the vertices of distance less than $r$ from the  unique cycle
of $M$, and (ii) contracting
$U$ into a single new vertex $u$. Suppose that $h'$ is a
homomorphism from $G'$ to $H$. Then by the previous paragraph,
$h'$ can be extended to a homomorphism from $G$ to $H$ where
each vertex $v\in V(G')-u$ is mapped to $h'(v)$, and every vertex in $U$
is mapped to $h'(u)$. Thus, if $G'$ is homomorphic to $H$ then so is $G$.

Let $T$ be
the tree which consists of a vertex adjacent to $L$
leaves. Since the degree of $u$ in $G'$ is at most $L$,
and using the fact that all vertices of $M$ are deleted
when forming $G'$ (here is where we require $r>|M|$),  the rest now
follows as in the $k\geq3$ case.
\proofend

\begin{lemma}
\label{distance} Let $H$ be a connected graph which contains a
triangle. Let $u$ be a vertex of $H$ and $M$ be a unicyclic graph
with unique cycle $C$. Denote the vertices of $M$ in a distance of
exactly $r \ge |V(H)|+3$ from $C$ by $U$. There is a homomorphism
from $M$ to $H$ such that all vertices in $U$ are mapped to $u$.
\end{lemma}
\begin{proof}
Let $h$ be a homomorphism from $C$ to the triangle $(v_1,v_2,v_3)$
of $H$. Observe that for $i=1,2,3$, there exist walks
$(v_i=)v_{i,0},\ldots,v_{i,r}(=u)$ of length exactly $r$ in $H$.
Let $w$ be a vertex in $M$ in the distance of $j \le r$ from $C$,
and $w'$ be the vertex of $C$ which has the distance $j$ from $w$.
Extend $h$ by assigning $h(w)=v_{ij}$ where $h(w')=v_i$. Observe
that $h$ is a partial homomorphism from $M$ to $H$ which maps
every vertex in $U$ to $u$. Trivially $h$ can be extended to a
homomorphism from $M$ to $H$.
\end{proof}

We close this section by discussing the possibilities of extending
Theorem \ref{thomo} to the case where $H$ is disconnected. We will
focus on graphs, i.e. the $k=2$ case. We can show that if
Hypothesis A does not hold for the $H$-homomorphism problem for
every graph $H$, then there must be a counterexample with two
components: a triangle and a graph $H_1$ that is triangle-free and
not 3-colourable. First note that every cycle is homomorphic to a
triangle, and a triangle is not homomorphic to any triangle-free
graph. So the condition of Hypothesis A is equivalent to  $H$
containing a triangle. On the other hand $H$ contains a
triangle-free component because being homomorphic to $H$ is
equivalent to being homomorphic to at least one of the components
$H_i$ of $H$, and so there is some component $H_i$, such that the
$H_i$-homomorphism problem has a coarse threshold. Let $H^1$ be
the subgraph of $H$ which consists of all triangles-free
components and $H^2$ be the remaining components of $H$. It is
easy to see that $H$ remains a counter-example if we add some
edges to $H^1$ without creating any triangle and we substitute
$H^2$ with a single triangle.

So the question of whether there is any graph $H$ for which the
$H$-homomorphism problem violates hypothesis A is equivalent to the following:

\begin{question}\label{q1} Is there any triangle-free graph $H$ with $\chi(H)>3$
such that for some values of $n$ and some $c>c(n)$, $G_{n,p=c/n}$
is not a.s. non-$H$-homomorphic, where $c(n)$ is the threshold value of
$3$-colorability?
\end{question}

\subsection{Directed Graphs}

Here, we provide an example of a directed graph $H$ for which
\begin{enumerate}
\item every unicyclic digraph has a homomorphism to $H$, and

\item the $H$-homomorphism problem under the $\csph$ model has a
coarse threshold.
\end{enumerate}

Unfortunately, this does not exhibit a coarse threshold under the
$\csp$ model, so the question of whether Hypothesis A holds for
all $H$-homomorphism problems for directed hypergraphs $H$ is
still open.

For a directed graph $D$, let $\tilde{D}$ denote the undirected
graph that is obtained from $D$ by ignoring the directions on the
edges. We define $D_{n,p}$ to be the random digraph on $n$
vertices where each of the $n(n-1)$ potential directed edges is
present with probability $p$.  Thus $D_{n,p}$ possible
constraining both $uv$ and $vu$ for some pair of vertices $v,u$,
i.e. a 2-cycle;  in fact, if $p=c/n$ for a constant $c$, then it
is straightforward to show that the probability that $D$ contains
at least one 2-cycle is $\z+o(1)$ for some constant $\z=\z(c)<1$.

$H$ consists of a specific digraph $H_1$, defined below, and a
pair of vertices $u_1,u_2$, where the edges $u_1,u_2$ and
$u_2,u_1$ are both present.
$H_1$ has the following properties:

\begin{itemize}
\item[(i):] every unicyclic digraph which does not contain a
2-cycle a.s. has a homomorphism to $H_1$ and

\item[(ii):] $D_{n,p=c_1/n}$ is not a.s. non-$H_1$-homomorphic for
some $c_1>1/2$.
\end{itemize}

It is easy to see that any unicyclic digraph, whose cycle is a
2-cycle, has a homomorphism to the 2-cycle.  By (i), every other
unicyclic digraph has a homomorphism to $H_1$.  Thus, every
unicyclic digraph has an $H$-homomorphism, as claimed. We will
show that, for every $1/2<c<c_1$, $D_{n,p=c/n}$ is neither a.s.
$H$-homomorphic nor a.s. non-$H$-homomorphic.  Thus, we have a
coarse threshold.  Condition (ii) above shows the latter, so we
just need to prove the former.

If $c>1/2$, then the graph $\tilde{D}$ for $D=D_{n,p}$, a.s. has a
giant component, as proven by Karp~\cite{karp}. It is not hard to
see that a.s. if $D$ has a 2-cycle in the giant component of
$\tilde{D}$, then there is no $H$-homomorphism: That 2-cycle must
be mapped onto $u_1,u_2$. Since $H$ has no edge incident to
$\{u_1,u_2\}$, any vertex that can be reached in $\tilde{D}$ from
that 2-cycle must also be mapped onto $u_1$ or $u_2$. So the
entire giant component must be mapped onto $\{u_1,u_2\}$. A.s.
that component has an odd cycle in $\tilde{D}$, and that odd cycle
cannot be mapped onto a 2-cycle. It is easy to show that the
probability that $\tilde{D}$ has a 2-cycle in its giant component
is at least some positive constant. Therefore, $D$ is not a.s.
$H$-homomorphic.

It remains only to prove the existence of some $H_1$ satisfying
(i), (ii). We choose $H_1$ to be a tournament (i.e. for every pair
of vertices, exactly one of the possible edges between them is
present) which contains every directed graph on $k_0$ vertices as
a subgraph where $k_0$ is a constant defined below.

For an undirected graph $G$ which does not contain any multiple
edges, the \emph{oriented chromatic number} $\chi_o$ of $G$ is the
minimum number $k$ such that every directed graph $D$ satisfying
$\tilde{D}=G$ is homomorphic to a directed graph $H$ with at most
$k$ vertices. The \emph{acyclic chromatic number} of a graph $G$
is the least integer $k$ for which there is a proper coloring of
the vertices of $G$ with $k$ colors in such a way that every cycle
of $G$ contains at least $3$ different colors. It was proved in
\cite{Raspaud} that if the acyclic chromatic number of a graph $G$
is bounded by $k$, then its oriented chromatic number is bounded
by $k.2^{k-1}$. When $D=D_{n,p}$ every edge is present in
$\tilde{D}$ with probability $2p-p^2$ and independent of the other
edges. Thus Lemma~\ref{acyclic} below together with the result
of~\cite{Raspaud} imply that taking $k_0=6 \times 2^{5}$ and
$c_1=c/2$, $H_1$ satisfies (ii), where $c$ is the constant which
is obtained from Lemma~\ref{acyclic}.

\begin{lemma}
\label{acyclic} There exists $c>1$ such that a.s. the acyclic
chromatic of $G_{n,p=c/n}$ is at most $5$.
\end{lemma}
\begin{proof}
Let $G=G_{n,p}$. A \emph{pendant} path in $G$ is a path in which
no vertices other than the endpoints lie in any edge of the graph
off the path. It is known that there exists $c>1$ such that a.s.
after removing the internal vertices of pendent paths of length at
least $4$ from $G$ every component is either a tree or it is
unicycle. One can use $3$ colors to color the vertices in these
components and then use $2$ other colors to color the removed
vertices such that every cycle in $G$ is colored by at least $3$
colors.
\end{proof}

\section{The $(d,k,t)$-model}\label{sdkt}


{\bf Proof of Theorem \ref{tdkt}} Suppose that the $(d,k,t)$-model
exhibits a coarse threshold. Then consider $p,\a,\e$ and $M$ as
guaranteed by Corollary~\ref{toolcsp}. It is easy to verify that,
since $t<d^{k-1}$ and $M$ is unicyclic, $M$ is satisfiable.
Suppose that $V(M)=u_1,...,u_r$ and let $a_i$ be the value of
$u_i$ in some particular satisfying assignment $A$ of $M$. Given a
CSP $F$ on at least $r$ variables, we define $F\oplus A$ to be the
CSP formed by choosing a random ordered $r$-tuple of variables
$v_1,...,v_r$ in $F$ and for each $1\leq i\leq r$, forcing $v_i$
to take the value $a_i$ by adding a one-variable constraint on
$v_i$. Clearly the probability that $\csph\oplus A$ is
unsatisfiable is at least as high as the probability that
$\csph\oplus M$ is unsatisfiable.


Lemma~\ref{local} implies that there exists some constant $L$ such that,
defining $E_1$ to be the event that every $v_i$ has at most $L$
neighbours, $\pr(E_1)\geq 1-\frac{\a}{2}$. Suppose that $E_1$
holds. Consider a particular $v_i$ and expose the $\ell\leq L$
edges containing it, $e_1,...,e_{\ell}$ and the corresponding
constraints $C_1,...,C_{\ell}$.  For each $C_j$, let $C'_j$ be the
$(k-1)$-variable constraint obtained by restricting $v_i$ to be
$a_i$; i.e. a $(k-1)$-tuple of values is permitted for $C'_j$ iff
$C_j$ permits that same $(k-1)$-tuple along with $v_i=a_i$. Since
$C_j$ has at most $t$ restrictions, $C'_j$ has at most $t$
restrictions.

Let $G$ be a random CSP formed as follows: start with a random
$\csph$ and then for each of the ${d^{k-1}\choose t}$ possible
constraints on $k-1$ variables and with $t$ restrictions, choose
$rL$ random ordered $(k-1)$-tuples of variables and place that
constraint on them. The probability that $G$ is unsatisfiable  is
at least as high as the probability that $\csph\oplus A$ is
unsatisfiable as each batch of $L$ copies of every
$(d,k-1,t)$-constraint is at least as restrictive as forcing
$v_i=a_i$. Thus, adding those $rL$ constraints boosts the
probability of unsatisfiability by at least $2\a$. We say that a
canonical set of ${d^{k-1}\choose t}rL$ ordered $(k-1)$-tuples is
{\em bad} if adding the constraints to that set results in an
unsatisfiable CSP. So, consider the following random experiment:
pick a random $\csph$ and then pick ${d^{k-1}\choose t}rL$ ordered
$(k-1)$-tuples of the variables. The probability that we pick a
bad set is at least $2\a$. Since ${d^{k-1}\choose t}rL=O(1)$, a
simple first moment calculation shows that a.s. the choice of
$(k-1)$-tuples will be vertex disjoint.  Thus, the probability of
picking a bad set is at least $2\a-o(1)$ even if we condition on
the $(k-1)$-tuples being vertex-disjoint.

Define $T$ by: (i) taking the hypergraph consisting of a vertex
$v$ lying in $rL{d^k\choose t}$ edges where no other vertex lies
in more than one of the edges, and (ii) placing each of the
${d^k\choose t}$ possible $(d,k,t)$-constraints on $rL$ of the
edges.  For each $1\leq\d\leq d$, let $T_{\d}$ denote the
collection of $(k-1)$-tuples obtained by removing $v$ from every
edge containing a constraint in which every restriction has
$v=\d$; note that $|T_{\d}|={d^{k-1}\choose t}rL$. By Lemma
\ref{tree}, there is some $\z=o(n^{1-k})$ such that
$\widehat{\CSP}_{n,p=z}({\cal P})$ a.s. contains a copy of $T$.
Consider adding that copy of $T$ to $\csph$. The probability that
for each $1\leq\d\leq d,$ $T_{\d}$ is a bad set is at least
$(2\a-o(1))^d$.  Note that if every $T_{\d}$ is a bad set, then
the resulting CSP is unsatisfiable because setting $v=\d$ requires
the set of $(k-1)$-constraints on $T_{\d}$ to be enforced. Thus,
the probability that $\widehat{\CSP}_{n,p=p+z}({\cal P})$ is
unsatisfiable is at least $(2\a-o(1))^d$.  By considering adding
$x$ copies of $T$, we see that the probability that
$\widehat{\CSP}_{n,p=p+xz}({\cal P})$ is satisfiable is at most
$(1-(2\a-o(1))^d)^x$ which is less than $\a$ for some sufficiently
large constant $x$.  Since $z=o(n^{1-k})$, this implies that
$\pr(\widehat{\CSP}_{n,(1+\e)p}({\cal P})\mbox{ is
unsatisfiable})>1-\a$ which contradicts
Corollary~\ref{toolcsp}(b). \proofend

\section{Binary CSP's with domain size $3$}\label{s23}

Recall that Istrate\cite{gi} (see also Creignou and
Daude\cite{cd2}) has proven that when the domain size $d=2$, then
Hypothesis A holds; i.e.  if every unicyclic CSP is satisfiable,
then $\csp$ has a sharp threshold.  This result does not extend to
$d=3$.  Consider the following example, with $d=3,k=2$:

\begin{example}\label{ed3}
We have two constraints.  $C_1$ says that either both variables are equal
to 1, or neither is equal to 1.  $C_2$ says that the variables cannot both
have the same value.
$\calp(C_1)=\frac{2}{3},\calp(C_2)=\frac{1}{3}$.
\end{example}

Observe that every unicyclic CSP that uses only constraints
$C_1,C_2$ is satisfiable.

Consider any $\frac{3}{2}<c<3$. Thus, a.s. the sub-CSP formed
by the $C_1$ constraints has a giant component,
and the sub-CSP formed by the $C_2$ constraints does not.
We will show that $\csp$
is neither a.s. satisfiable nor a.s. unsatisfiable.

To see that it is not a.s. unsatisfiable, note that the subgraph
induced by the $C_2$ constraints is 2-colourable with probability
at least some positive constant. This follows from the well known
fact that the for $c<1$ the random graph $G(n,\frac{c}{n})$ is
2-colourable with probability at least some positive constant. If
it is 2-colourable, then we can satisfy all the $C_2$ constraints
by assigning every variable either 2 or 3; this will not violate
any $C_1$ constraints.

To see that it is not a.s. satisfiable, note that subgraph formed
by the $C_1$ constraints has a giant component $T$. So either
every variable in $T$ is assigned 1 or none of them are. A.s. at
least one $C_2$ constraint has both variables in $T$, and so they
cannot both be assigned 1. Thus, a.s. no variables in $T$ can be
assigned 1. This implies that if the $C_2$ constraints form an odd
cycle using variables of $T$ then the CSP is not satisfiable; that
event occurs with probability at least some positive constant.

The main result of this section, is that when $d=3$ and $k=2$,
if Hypothesis A fails on some model, then it has to fail
for the same reason as it failed for Example \ref{ed3}.

Consider a CSP $F$ where every constraint is on 2 variables.
Suppose there is some constraint on variables $v,u$ which
implies that if $v$ is assigned $\d$ then $u$ must be assigned $\g$;
we say that {\em $v:\d$ forces $u:\g$} and denote this by
$v:i\rightarrow u:j$. Moreover if there is a sequence
of variables $v_1,\ldots,v_r$ and values
$\d_1,\ldots,\d_r$ such that $v_i:\d_i\rightarrow v_{i+1}:\d_{i+1}$
for $i=1,\ldots,r-1$ then we say that $v_1:\d_1$ forces $v_r:\d_r$.

\begin{theorem}\label{t23} Consider some $\calp$ with
$d=2,k=3$ such that every unicyclic CSP formed from
$\supp(\calp)$ is satisfiable.  $\csp$ has a coarse threshold
iff there exists a unicyclic CSP $M$
formed from $\supp(\calp)$,
a value $1\leq\d\leq3$, $p=p(n),\e>0,z>0,b>0$ such that:
\begin{enumerate}
\item[(a)] $\e<\pr(\csp\mbox{ is satisfiable })<1-\e$;

\item[(b)] $M$ cannot be satisfied using only the two values other
than $\d$;

\item[(c)] $\csph$ a.s. has at least $z n$ variables $v$ such that
$v:\d\rightarrow u:\d$ for at least $b n$ variables $u$.
\end{enumerate}
\end{theorem}

Thus, this explains the only ways that a model with $d=3,k=2$
can have a coarse threshold.  We remark that this theorem does not
extend to $d=4,k=2$ nor $d=3,k=3$; in both cases there are other causes
for a coarse threshold.

{\bf Proof}
We leave it to the reader to verify, using similar reasoning to
that for Example \ref{ed3} that conditions (a,b,c) imply a coarse threshold.
The only difference here is that the set of variables $u$
in (c) could change for different choices of $v$, and it is important to
note that for each value $\d$, there is some constraint in $\supp(\calp)$
that forbids both variables from receiving $\d$ as otherwise
$\csp$ is trivially satisfiable by setting every variable equal to $\d$.
We will give a proof of the other direction.
The case where some domain values are {\em bad} (as defined in \cite{mm1})
is easily disposed of, so we assume that there are no such values.

For each variable $v$ and each $1\leq \d\leq3, 1\leq \g\leq 3$
we define $F_{\d,\g}(v)$ to be the set of variables $u$ such that
$v:\d\rightarrow u:\g$, and we define $F_{\d}(v)=\cup_{1\leq\g\leq3}F_{\d,\g}(v)$.
We can expose $F_{\d}(v)$ by using a simple breadth-first
search from $v$.  This allows us to analyze the distribution of the
size of $F_{\d}(v)$ and $F_{\d,\g}(v)$ using a standard branching-process analysis
(see eg.\ Chapter 5 of \cite{jlr}). We say that $F_{\d,\g}$ {\em percolates}
if there are constants $\z,\b>0$ such that
$\pr(|F_{\d,\g}(v)|\geq \b n)\geq \z$. It is straightforward to prove:

{\em Claim 1: If $F_{\d,\g}$ does not percolate, then for every $\xi>0$
there is a constant $L$ such that\\ $\pr(|F_{\d,\g}(v)|\leq L)>1-\xi$.}

{\em Claim 2: If $F_{\d,\g}$ percolates, then there are constants $z,b>0$
such that a.s. there are at least $z n$ variables $v$ with $|F_{\d,\g}(v)|\geq\b n$.}

Claim 2 yields that $F_{\d,\d}$ percolates for at most one value
$\d$:  Suppose that $F_{\d,\d}$ and $F_{\g,\g}$ both percolate. We
want to show that in this case Corollary~\ref{toolcsp}(b) fails.
To this end, we obtain an instance of
$\widehat{\CSP}_{n,p(1+\epsilon)}({\cal P})$ as follows. First, we
consider $P_0$, an instance of $\widehat{\CSP}_{n,p}({\cal P})$.
 Then $P_i$ is obtained by taking the union of $P_{i-1}$ and an
instance of  $\widehat{\CSP}_{n,\frac{p\epsilon}{3}}({\cal P})$,
for $i=1,2,3$.

 Let $u$, $v$, and $w$ be variables such that in $P_0$, $u:\d$
forces both $v:\d$ and $w:\d$. If the restriction $(\d,\d)$ is
added on $v$ and $w$, then $u$ cannot be assigned $\d$. Since
there are constraints in $\supp({\cal P})$ with the restriction
$(\d,\d)$ and constraints with the restriction $(\g,\g)$ (see
\cite{mm1}), we can conclude that in $P_1$, a.s. there are two
sets of variables of size $\Theta(n)$, $A$ and $B$ such that $\d$
cannot be assigned to any variable in $A$, and $\g$ cannot be
assigned to any variable in $B$.

Since $F_{\d,\d}$ percolates, there is a value $\sigma$ and a
constraint in $\supp({\cal P})$
such that if that constraint is applied on $v_1,v_2$, then
$v_1:\sigma \rightarrow v_2:\d$. Let $u$ be a variable in $A$,
and $v$ be an arbitrary variable. If there is a constraint on $u$
and $v$ which justifies that $v:\sigma \rightarrow u:\d$, then
$\sigma$ cannot be assigned to $v$. Now we can conclude that in
$P_2$, almost surely, there is a set $C$ of size $\Theta(n)$ of
variables such that $\sigma$ cannot be assigned to any variable in
$C$ for the reason mentioned above. Every variable in $A$ or $B$ is in $C$ with a positive
probability and independent of the other variables. So $C_1=C\cap
A$ and $C_2=C \cap B$ are both of size $\Theta(n)$, almost surely.
Without loss of generality suppose that $\delta \neq \sigma$
(otherwise we would assume $\g \neq \sigma$). Both of the values
$\delta$ and $\sigma$ cannot be assigned to any variables in
$C_1$. So there is only one value which can be assigned to these
variables. But then in $P_3$ a.s. there is a constraint on two
variables in $C_1$ which forbids them to be assigned this value
simultaneously. So $\widehat{\CSP}_{n,p(1+\epsilon)}({\cal P})$ is
not satisfiable a.s. which contradicts Corollary~\ref{toolcsp}(b).

It is less straightforward, but not very difficult, to prove:

{\em Claim 3: If $F_{\d,\g}$ percolates for any pair $\d,\g$ then
either (i) $F_{\d,\d}$ percolates or (ii) there is some $\m$ such that $F_{\m,\m}$
percolates and there is a sequence of constraints in $\supp(\calp)$
through which $v:\d\rightarrow u:\m$.}

Suppose that $\csph$ has a coarse threshold and consider
$M,\e,\a,p=p(n)$ from Corollary~\ref{toolcsp}.

{\em Claim 4: There is a value $\d$ such that (i)
that every satisfying assignment of $M$ must use $\d$ on at least one variable and (ii)
$F_{\d,\g}$ percolates for at least one value $\g$.}

If $F_{\d,\d}$ percolates
then this satisfies Theorem \ref{t23}.  Otherwise, applying Claim 3, $F_{\m,\m}$
percolates and there is a sequence of constraints so that $v:\d\rightarrow u:\m$.
Attaching that sequence to every variable of $M$ yields a unicyclic CSP $M'$
for which every satisfying assignment must use $\m$ on at least one variable.
Thus  $M',\m$ satisfy Theorem \ref{t23}.

Suppose that Claim 4 does not hold, and consider any satisfying assignment $A$ of $M$
in which every value $\d$ used is such that $F_{\d,\d}$ does not percolate.  Suppose
that $M$ has $r$ variables $x_1,...,x_i$ and that $A$ assigns $a_i$ to $x_i$.
Recall from Section \ref{sdkt} that $\csp\oplus A$ is formed by taking
$\csp$ and then choosing $r$ random variables $v_1,...,v_i$ and adding
one-variable constraints that force $v_i$ to take $a_i$.  Clearly
$\pr(\csp\oplus A\mbox{ is unsatisfiable})\geq\pr(\csp\oplus M\mbox{ is unsatisfiable})$.

Expose $F=\cup_{i=1}^rF_{a_i}(v_i)$, and $U$, the set of variables
outside of $F$ that lie in a constraint with a variable in $F$.
Since none of the $F_{a_i}$ percolate, Claim 2
allows us to show that there is some $L$ such that with probability at least $1-\a/2$,
$|U|<L$.  Since adding $M$ to $\csph$ increases
the probability of unsatisfiability by at least $2\a$, it must be that the probability
that $\csph$ is satisfiable, $|U|\leq L$ and $\csph\oplus A$ is unsatisfiable
is at least $3\a/2$.

Suppose that $\csph$ is satisfiable and $|U|\leq L$.   Consider some $u\in U$
sharing a constraint with $w\in F$ where $A$ forces $w$ to take the value $\m$.
Let $\Omega=\Omega(u)$ be the set of values which can be assigned to $u$ which,
in conjunction with assigning $\m$ to $w$ do not violate their constraint.
We know that $|\Omega|\neq0$ since otherwise $\m$ is a {\em bad value}.
We know that $|\Omega|\neq1$ since otherwise $u\in F$. So
$|\Omega(u)|\geq 2$ for each $u\in U$. Suppose that $u_1,...,u_\ell$
are the variables in $U$ with $|\Omega|=2$,
and let $\d_i$ be the value not in $\Omega(u_i)$. Consider taking
a random CSP formed as follows:  first take a $\csph$ and then
choose $\ell$ random variables $u_1,...,u_{\ell}$ and force $u_i$ to not
take value $\d_i$ using a one-variable constraint.
We have proved is that the one-variable constraints
boost the probability of unsatisfiability by $3\a/2$.

At this point, we can complete the proof using the argument from
the Achlioptas-Friedgut proof\cite{af} that $d$-colourability has a sharp threshold
for $d\geq 3$.  In that paper, the starting point was to note
that since every constraint is a colourability constraint,
fixing an assignment on $M$ at worst forbids one colour from
each neighbour of $M$.  This put them in the same position that
we are in now, and so the rest of our proof is the same as theirs.
\proofend

We close this section by noting why this proof can not be extended
to general $d$.  The problem is that possibly some of the
variables in $U$ would have their domain sizes reduced by two
instead of one. The argument in \cite{af} cannot handle that
possibility.

\section{Future Directions}
There is clearly much work still to be done along these lines of research.
The big problem still remains - determine precisely
which models from \cite{mm1} have a sharp threshold.
Of course, Section 2 indicates that this may be overly ambitious.
But lowering our sights only slightly, we can try to determine
all possible {\em causes} for coarse thresholds, i.e. continue the
course started in Section \ref{s23}.  An important subgoal would
be to do this for binary CSP's, i.e. the case where $k=2$.
Another reasonable goal to pursue would be to cover the $d=3$ case.

As far as more specific classes of models go, one should try to
extend the work in Section \ref{sh} and examine
whether Hypothesis A holds for $H$-homomorphism problems
when $H$ is a {\em directed} hypergraph.  Such homomorphism
problems are equivalent to CSP's in which every constraint
is identical under some permutation of the variables.
And of course,
it would be good to determine whether the ``connected'' condition can
be removed from Theorem \ref{thomo}.

\end{document}